\documentstyle{article}
\date{}
\title{New  Existence Theorems about the Solutions of Some Stochastic Integral Equations}
\author{{Xuemei Chen  and Yingying Qi and Chunyan Yang}\\
{\small Department of Mathematics, Sichuan University,}\\
{\small Chengdu 610064, P.R.China}}
\begin{document}
\maketitle \large \baselineskip 14pt

\begin{quote}

{\bf Abstract}\ \ Picard's iteration has been used to prove the
existence and uniqueness of the solution for stochastic integral
equations, here we use Schauder's fixed point theorem to give a new
existence theorem about the solution of a stochastic integral
equation, our theorem can weak some conditions gotten by applying
Banach's fixed point theorem.

{\bf Keywords}\ \ Stochastic integral equation; Schauder's fixed
point theorem; bounded closed convex subset; compact operator

{\bf AMS Subject Classification}\ \ 47H10, 60H05, 60H10

\end{quote}

\section{Introduction and Main Results}

\ \ \ \ \  Many mathematical models of phenomena occuring in
sociology, physics, biology and engineering involve random
differential and integral equations. Theoretical and applied
treatments of problems concerning random differential and integral
equations can be found in many papers and monographs: Bharucha-Reid
([4]); Doob ([5]); Padgeett and Tsokos ([6]); Tsokos and Padgett
([7]); Rao and Tsokos ([8]).
\\

For stochastic differential equations, firstly we should prove the
existence for the solutions. In filtering problem we know the system
$X_t$ satisfying $dX_t=b(t,X_t)dt+\sigma (t,X_t)dB(t),$ the
observations $Z_s$ satisfying $dZ_s=c(s,X_s)ds+d(s,X_s)dV(s),Z_0=0 $
for $0\leq s\leq t$, where V(s) is Brownian motion, we want to find
the best estimate $\hat{X} _t$ of the state $X_t$ of the system
based on these observations, before we find the estimate $\hat{X}
_t$, we should give some assumptions for
 the existence of the corresponding stochastic integral equations.
\\

Some mathematicians have used Picard's iteration or Banach's
contraction mapping principle to prove the existence and uniqueness
of some stochastic integral equations. The goal of this paper is to
give a new existence theorem about stochastic integral equations
using Schauder's fixed point theorem, in order to apply Schauder's
fixed point theorem we need to
 construct a compact operator A and a convex bounded and closed
nonempty subset M. Furthermore, comparing with Banach's fixed point
theorem we weak some conditions.

\vspace{0.3cm} \textbf{Theorem 1.1(\cite{1})}(Schauder's fixed point
theorem).\ {\em The compact operator $$A: M\longrightarrow M$$ has
at least one fixed point when M is a bounded, closed, convex,
nonempty subset of a Banach space X over real field.}

\vspace{0.3cm} \textbf{Theorem 1.2(\cite{1})}(Banach's fixed point
theorem).{\em We assume that:\\

(a) M is a closed nonempty subset in the Banach space X over field K, and\\

(b) the operator $A: M\longrightarrow M$ is k-contractive, i.e.,
there is $0\leq k<1$ such that:
$$
{\| Au-Av\| \leq k\| u-v\|   \ \ for \ all \ u,v\in M.}
$$

Then, the operator A has exactly one fixed point $u$ on the set M.}\\

Schauder's Fixed Point Theorem can be applied to many fields in
mathematics, especially to the integral equation:
\begin{eqnarray}
{u(x)=\lambda\int_a^bF(x,y,u(y))dy,  \ \ \ a\leq x\leq
b,\label{91.1}}
\end{eqnarray}
where $-\infty <a<b<+\infty$ and $\lambda\in R$. Let
$$
Q={\{(x,y,u)\in R^3 ;\ x,y\in [a,b],\ |u|\leq r\}} \ \ for \ fixed \
r>0.$$

\vspace{0.3cm} \textbf{Theorem 1.3(\cite{1})}{\em  Assume the
following conditions:\\

(a) The function $F:Q\longrightarrow R $ is continuous.\\

(b) We define $(b-a)$$\mathcal {M}$ $ : =\max  _{(x,y,u)\in Q}$ \
$|F(x,y,u)|$, let the real number $\lambda $ be given such that
\begin{eqnarray}
|\lambda | {\mathcal {M}}\leq r
\end{eqnarray}

(c) We set $X: =C[a,b]$ and $ M:=\{u\in X;\| u\| \leq
r\}$.\\

Then the original integral equation (1) has at least one solution
$u\in M$.}\\

 It well known that Banach's Fixed Point Theorem can be used to prove
the existence and uniqueness of the solution for the integral
equation(1).

\vspace{0.3cm} \textbf{Theorem 1.4(\cite{1})}{\em Assume the
following conditions:\\

(a) The function $F:[a,b]\times [a,b]\times R\longrightarrow R$ is
continuous, and the partial derivative
$$F_u:[a,b]\times [a,b]\times R\longrightarrow R$$ is also
continuous. \\

(b) There is a number $\mathcal {L}$ such that
$$
|F_u(x,y,u)|\leq {\mathcal {L}} \ \ for \ all \ x,y\in [a,b],u\in R.
$$

(c) Let the real number $\lambda $ be given such that
\begin{eqnarray}
(b-a)|\lambda |{\mathcal {L}}<1.
\end{eqnarray}

(d) Set $X:=C[a,b]$ and $\| u\| :=max_{a\leq x\leq
b}|u(x)|$.\\

Then the original equation (1) has a unique solution $u\in M$.}\\

From Theorem 1.1 and Theorem 1.2 we know that Schauder's fixed point
theorem is a existence principle, while Banach's fixed point theorem
is a existence and uniqueness theorem. It seems that the conditions
in Theorem 1.3 are weaker comparing with Theorem 1.4, that is, (2) \
is \ easier \ to \ reach \ than \ (3).\\

Here, a nature question is whether Schauder's and Banach's fixed
point theorems can be applied to stochastic integral equations.
Furthermore, whether the conditions coming from Banach's
fixed point theorem are stronger than the conditions from Schauder's fixed point theorem.\\

\vspace{0.3cm} {\bf Notation 2.1(\cite{2})}\ {\em For convenience,
we will use $X:=L_{ad}^2([a,b]\times \Omega )$ to denote the space
of all stochastic process $f(t,\omega ), a\leq t\leq b,\omega \in
\Omega
$ satisfying the following conditions:\\

(1)$f(t,\omega )$ is adapted to the filtration $\mathcal {F}$$_t$; \\

(2)$\int _a^b E(|f(t)|^2)dt<+\infty .$}
\\

Now we want to solve the stochastic integral equation:
\begin{eqnarray}
{x(t;w)=h(t_0;w)+\int_a^t \sigma (s,x(s;w))dB(s)+\int_a^t
f(s,x(s;w))ds,\  a\leq t\leq b,}
\end{eqnarray}
where:

(1)$\omega \in \Omega $ , where $\Omega $ is the supporting set of
the probability measure space $(\Omega ,{\mathcal {F}} ,P)$ with
$\mathcal {F} $
being the $\sigma $ -algebra and P the probability measure;\\

(2)$x(t;w)$ is the unknown random variable for each $t\in [a,b]$;\\

(3)$h(t_0;w)$ is the knowm random variable and $E|h|^2<+\infty $;\\

(4)$B(t)$ be a Brownian motion and $\{{\mathcal {F}} _t ; a\leq
t\leq b\}$ be a filtration so there $B(t)$ is ${\mathcal {F}} _t$
-measurable for each $t$ and $B(t)-B(s)$ is independent of
${\mathcal {F}}_s$ for any $s<t$.

Let
$$
Q={\{(t,X_t)\in R^2 ;\ t\in [a,b], \ and  \ \| X_t\| \leq r  \ for \
fixed
 \ r>0\}}.
$$

\vspace{0.3cm} \textbf{Theorem 1.5 }{\em Assume the
following conditions: \\

(a) $f(s,X_s),\ \sigma (s,X_s)$ are measurable on $[a,b]\times
\Omega
$;\\

   ~~~~ $f(s,X_s):Q\longrightarrow R$ is continuous;\\

   ~~~~ $\sigma (s,X_s):Q\longrightarrow R$ is continuous;\\

(b) We define
\begin{eqnarray}
d=\sup _{(s,X_s)\in Q}\{ \| f(s,X_s)\| ,\| \sigma (s,X_s)\| \},
\end{eqnarray}
 and let the real number $a,\ b,\ d $ and
random variable $h$ be given that
\begin{eqnarray}
3E[h^2]+3(1+b-a)(b-a)d^2\leq r^2.
\end{eqnarray}

(c) We set $X:=L_{ad}^2([a,b]\times \Omega )$ (see Notation 2.1) and
$M:=\{ X_t\in X;\| X_t\| \leq r\}$.\\

Then the stochastic integral equation (4) has at least one solution
$X_t\in M$}.

\vspace{0.3cm} \textbf{Theorem 1.6 }{\em Assume the
following conditions: \\

(a) $f(s,X_s), \ \sigma (s,X_s)$ are measurable on $[a,b]\times
\Omega
$;\\

(b)\begin{eqnarray}
 \  |f(s,X_s)-f(s,Y_s)|&\leq& k_1 |X_s-Y_s|;\\
|\sigma (s,X_s)-\sigma (s,Y_s)|&\leq& k_2 |X_s-Y_s|;
 \end{eqnarray}

(c) Let the real number $a,\ b,\ c=\{k_1,k_2\}$ be given such that
\begin{eqnarray}
0\leq 2c^2(1+b-a)(b-a)<1.
\end{eqnarray}

(d) We set $X:=L_{ad}^2([a,b]\times \Omega )$ (see Notation 2.1) and
$M:=\{ X_t\in X;\| X_t\| \leq r\}.$ \\

Then the stochastic integral equation (4) has a unique solution
$X_t\in M$}.

\section{Some Lemmas }
\ \ \ \ \ \ We require the following Lemmas for proving the
existence of the stochastic integral equation.

\vspace{0.3cm} {\bf Lemma 2.1(\cite{1})}\ {\em Let X and Y be normed
spaces over field K, and let $$A:M\subseteq X\longrightarrow Y$$ be
a continuous operator on the compact nonempty subset M of X. Then, A
is uniformly continuous on M.

\vspace{0.3cm} {\bf Lemma 2.2(\cite{1})}\  Let
$X:=L_{ad}^2([a,b]\times \Omega )$ with $\|
X_t\|:=(E|X_t|^2)^\frac{1}{2}$ and $-\infty
<a<b<+\infty $. Suppose that we are given a set M in X such that \\

(1)M is bounded , i.e., $\| X_t\| \leq r$ for all $X_t\in M$ and
fixed $r\geq 0$.\\

(2)M is equicontinuous, i.e., for each $\varepsilon
>0$, there is a $\delta >0$ such that \\
$$|t_1-t_2|<\delta \ \ \ and \ \ X_t\in M \ \ \ imply
\ \ |X_{t_1}-X_{t_2}|<\varepsilon .$$ Then, M is a relatively
compact subset of X.

 \vspace{0.3cm} {\bf Definition 2.1(\cite{1})}\
Let X and Y be normed spaces over field K. The operator
$$A:M\subseteq
X\longrightarrow Y$$ is called compact iff\\

(1) A is continuous, and\\

(2) A transforms bounded sets into relatively compact sets.

\vspace{0.3cm} {\bf Lemma 2.3(\cite{3})}\ (It$\hat{o}$ Isometry)\\

For each $X_t,Y_t\in L_{ad}^2([a,b]\times \Omega )$, we have\\
$$E[(\int _a^b f(t,w)dB(t))^2]=E[\int _a^b f^2(t,w)dt]$$}

\section{The Proof of  Theorem 1.5}
\ \ \ \ \ \ \vspace{0.3cm}  {\bf Proof}: We divide the proof into
three steps:

$Step$ 1: We prove that $M=\{ X_{t}\in X;\| X_t\| \leq r\}$ is
closed, convex subset of $L_{ad}^2([a,b]\times \Omega )$ (see
Notation 2.1).
\\

(A) We prove M is closed. \\

Let $X_t^{(n)}\in M$ for all $n$, i.e., $$\| X_t^{(n)}\| \leq r \ \
\ for  \ all \  n.$$ If $X_t^{(n)}\longrightarrow X_t$ as
$n\longrightarrow +\infty $, then
$\|X_t\|\leq r$, and hence $X_t\in M$.\\

(B) We prove M is convex. \\

If $X_t,Y_t\in M$ and $0\leq \alpha \leq 1$, then
\begin{eqnarray*}
\ \|\alpha X_t+(1-\alpha )Y_t\| &\leq&\| \alpha X_t\|
+\|(1-\alpha ) Y_t\| \\
&\leq& \mbox{} \alpha r+(1-\alpha )r\\
&=& \mbox{} r
\end{eqnarray*}
Hence $$\alpha X_t+(1-\alpha )Y_t\in M.$$

$Step$ 2:  We prove that $A: M\longrightarrow M$ is a compact
operator.\\

Define $A: M\longrightarrow M$
\begin{eqnarray}
{A(X_t)=h(t_0;w)+\int_a^t \sigma (s,x(s;w))dB(s)+\int_a^t
f(s,x(s;w))ds,  \ a\leq t\leq b,}
\end{eqnarray}
Then\\

(a) $A: M\longrightarrow M$ is a continuous operator.\\

 By Lemma
2.1, we know $f(s,X_s), \sigma (s,X_s)$ are uniformly continuous on
the compact set Q. This implies that, for each $\varepsilon
>0$, there is a number $\delta
>0$ such that
$$\| \sigma (s,X_s)-\sigma (s,Y_s)\| <\varepsilon _1$$
$$\| f(s,X_s)-f(s,Y_s)\| <\varepsilon _2$$
for all$(s,X_s),(s,Y_s)\in Q$ with $\| X_s-Y_s\|
<\delta .$\\

For each $X_t,Y_t\in M$, we have
\begin{eqnarray*}
\ \| AX_t-AY_t \| ^2&=&E(\int _a^t (\sigma
(s,X_s)-\sigma (s,Y_s))dB(s)\\
&+& \mbox{} \int _a^t (f(s,X_s)-f(s,Y_s))ds)^2
\end{eqnarray*}
Using the inequality $(a+b)^2\leq 2(a^2+b^2)$ to get
\begin{eqnarray}
\ \| AX_t-AY_t \| ^2&\leq&2E(\int _a^t (\sigma
(s,X_s)-\sigma (s,Y_s))dB(s))^2{\nonumber}\\
&+& \mbox{} 2E(\int _a^t (f(s,X_s)-f(s,Y_s))ds)^2\label{95}
\end{eqnarray}
Applying the It$\hat{o} $ Isometry to $E(\int _a^t (\sigma
(s,X_s)-\sigma (s,Y_s))dB(s))^2$ , we get:
\begin{eqnarray}
\ E(\int _a^t (\sigma (s,X_s)-\sigma (s,Y_s))dB(s))^2&=&E(\int _a^t
(\sigma (s,X_s)-\sigma (s,Y_s))^2ds){\nonumber}\\
&=& \mbox{} \int _a^t E(\sigma (s,X_s)-\sigma (s,Y_s))^2ds{\nonumber}\\
&=& \mbox{} \int _a^t \| \sigma (s,X_s)-\sigma (s,Y_s)\|
^2ds{\nonumber}\\
&<&  \mbox{} (b-a)\varepsilon _{1}^2\label{95}
\end{eqnarray}
For $E(\int _a^t (f(s,X_s)-f(s,Y_s))ds)^2$, we use Schwarz's
inequality to get
\begin{eqnarray}
\ E(\int _a^t (f(s,X_s)-f(s,Y_s))ds)^2&\leq&E((t-a)\int _a^t
(f(s,X_s)-f(s,Y_s))^2ds){\nonumber}\\
&\leq& \mbox{} (b-a)\int _a^t E(f(s,X_s)-f(s,Y_s))^2ds{\nonumber}\\
&=& \mbox{} (b-a)\int _a^t \| f(s,X_s)-f(s,Y_s)\|
 ^2ds{\nonumber}\\
 &<& \mbox{} (b-a)^2\varepsilon _{2}^2\label{95}
\end{eqnarray}
Put equations (12) and (13) into equation (11) to get
\begin{eqnarray*}
\ \| AX_t-AY_t \| ^2&<&2(b-a)[(b-a)\varepsilon _{2}^2+\varepsilon _{1}^2]\\
&=& \mbox{} \varepsilon ^2
\end{eqnarray*}

Therefore for each $X_t,Y_t\in M$, there exists $\delta >0$, when
$\| X_t-Y_t\| \leq \delta $, we have
\begin{eqnarray*}
\| AX_t-AY_t \| < \varepsilon .
\end{eqnarray*}
That is: A is a continuous operator. \\

(b) A(M) is bounded.\\

 For each
$X_t\in M$
\begin{eqnarray}
\ \| AX_t\| ^2&=&E(h(t_0;\omega )+\int _a^t \sigma (s,X_s)dB(s)+\int _a^t f(s,X_s)ds)^2{\nonumber}\\
 &\leq& \mbox{} 3E[h^2]+3E(\int _a^t \sigma
(s,X_s)dB(s))^2+3E(\int _a^t f(s,X_s)ds)^2{\nonumber}\\
 &\leq& \mbox{} 3E[h^2]+3\int _a^t E|\sigma
(s,X_s)|^2ds+3(b-a)\int _a^t E|f(s,X_s)|^2ds{\nonumber}\\
&=& \mbox{} 3E[h^2]+3\int _a^t \| \sigma
(s,X_s)\| ^2ds+3(b-a)\int _a^t \| f(s,X_s)\| ^2ds{\nonumber}\\
&\leq& \mbox{} 3E[h^2]+3(b-a)(1+b-a)d^2{\nonumber}\\
&\leq& \mbox{} r^2
\end{eqnarray}
Thus A(M) is bounded.\\

(c) A(M) is equicontinuous. \\

For each $X_t\in M$, we have
\begin{eqnarray*}
\ \| AX_{t_1}-AX_{t_2} \| ^2&=&E(\int _{t_2}^{t_1} \sigma
(s,X_s)dB(s)+\int _{t_2}^{t_1} f(s,X_s)ds)^2\\
&\leq& \mbox{} 2E(\int _{t_2}^{t_1} \sigma (s,X_s)dB(s))^2+2E(\int
_{t_2}^{t_1} f(s,X_s)ds)^2\\
&\leq& \mbox{} 2\int _{t_2}^{t_1} E|\sigma (s,X_s)|^2ds+2(b-a)\int
_{t_2}^{t_1} E|f(s,X_s)|^2ds\\
&\leq& \mbox{} 2(1+b-a)|t_1-t_2|d^2\\
\end{eqnarray*}
Take $$\delta =\frac{\varepsilon ^2}{2(1+b-a)d^2},$$ then for each
$\varepsilon
>0$, there exists $$\delta =\frac{\varepsilon ^2}{2(1+b-a)d^2},$$ when
$|t_1-t_2|<\delta $, we have
$$\|
AX_{t_1}-AX_{t_2}\| < \varepsilon.$$
Hence A(M) is equicontinuous.\\

Then by Lemma 2.2 and Definition 2.1, we know $A:M\longrightarrow M$
is a compact
operator.\\

$Step$ 3:  we prove that $A(M)\subseteq M$.\\

For each $X_t\in M$, we have
\begin{eqnarray*}
\ \| AX_t\| ^2&=&E(h+\int _a^t \sigma
(s,X_s)dB(s)+\int _a^t f(s,X_s)ds)^2\\
&\leq& \mbox{} 3E[h^2]+3E(\int _a^t \sigma
(s,X_s)dB(s))^2+3E(\int _a^t f(s,X_s)ds)^2\\
 &\leq& \mbox{} 3E[h^2]+3\int _a^t E|\sigma
(s,X_s)|^2ds+3(b-a)\int _a^t E|f(s,X_s)|^2ds
\end{eqnarray*}
where $f(s,X_s),\sigma (s,X_s)\in L_{ad}^2([a,b]\times \Omega )$.

 So $$\int _a^b
\| AX_t\| ^2dt<+\infty ,$$\\
that is $$ AX_t\in L_{ad}^2([a,b]\times \Omega ),$$ meanwhile, we
have proved $\| AX_t\| \leq r$  in (14), therefore $AX_t\in M$, that
is $A(M)\subseteq M$.
\\

Thus the Schauder's Fixed Point Theorem tells us that equation (4) has at least one solution $X_t\in M$.\\

\section{The Proof of the Theorem 1.6}
\ \ \ \ \ \ \vspace{0.3cm}  {\bf Proof}: \ In the proof of Theorem
1.5, we have proved that M is
closed.\\
We now show that A is a contractive mapping:\\

For each $X_t,Y_t\in M$, we have
\begin{eqnarray*}
\ \| AX_t-AY_t \| ^2&=&E(\int _a^t (\sigma (s,X_s)-\sigma
(s,Y_s))dB(s) + \mbox{} \int _a^t (f(s,X_s)-f(s,Y_s))ds)^2\\ &\leq&
\mbox{} 2E(\int_a^t (\sigma (s,X_s)-\sigma (s,Y_s))dB(s))^2
+2E(\int _a^t (f(s,X_s)-f(s,Y_s))ds)^2\\
&\leq& \mbox{} 2\int _a^t E(\sigma (s,X_s)-\sigma
(s,Y_s))^2ds+2(b-a)\int _a^t E(f(s,X_s)-f(s,Y_s))^2ds\\
 &\leq& \mbox{} 2k_2^2\int _a^t
E|X_s-Y_s|^2ds+2k_1^2(b-a)\int _a^t
E|X_s-Y_s|^2ds\\
&\leq& \mbox{} 2c^2(1+b-a)(b-a)\| X_t-Y_t\| ^2
\end{eqnarray*}
Let $$k^2=2c^2(1+b-a)(b-a)<1,$$ then\\ $$\| AX_t-AY_t\| \leq k\|
X_t-Y_t\|, \ \ 0\leq k<1 ,$$
Therefore A is k-contractive.\\

Then the Banach's fixed point theorem tells us that the stochastic
integral equation (4) has a unique solution $X_t\in M$.

\section{Comparing Theorem 1.5 with Theorem 1.6}
\ \ \ \ \ \ Comparing Theorem 1.5 with Theorem 1.6 we know when we
use Schauder's fixed point theorem to prove the existence of the
solution for the integral equation, we need conditions (5) and (6).

But when we use Banach's fixed point theorem to prove the existence
of the solution for the stochastic integral equation, we need
conditions (7), (8) and (9).

Obviously, the condition (6) is weaker than the condition (9).

\section{An Example for Theorem 1.5}
 \ \ \ \ \ \ We apply the above Theorem 1.5 to the linear stochastic
 integral equation:
\begin{eqnarray}
X_t=\int _0^t f(s)X_sds+\int _0^t g(s)X_sdB(s), \ \ 0\leq t\leq 1,
\end{eqnarray}

\vspace{0.3cm}  {\bf Proof}: Define the operator:
   $$A(X_t)=\int _0^t f(s)X_sds+\int _0^t g(s)X_sdB(s), \ \ 0\leq t\leq 1,$$

Obviously, $f(s)X_s$ and $g(s)X_s$ are measurable and continuous. We
define that $d=\sup E|X_s|^2$, $c=\max \{f(s)^2,g(s) ^2\}$,
$6cd=r^2$ and set $X:=L_{ad}^2([a,b]\times
\Omega )$ and $M=\{ X_t\in X;\| X_t\| \leq r\}$.\\

Then all conditions of Theorem 1.5 hold, so (15) has at least one
solution.\\

Especially, if we let $f(s)=u,g(s)=\sigma $, then the equation is
the geometric Brownian motion equation
\begin{eqnarray}
X_t=\int _0^t uX_sds+\int _0^t \sigma X_sdB(s), \ \ 0\leq t\leq 1,
\end{eqnarray}
where $u$ is the expected return rate(constant), $\sigma $ is
volatility(constant), $B(t)$ is standard Brown motion.\\

We can easily to get that the equation (16) has at least one
solution. It well known that the existence of the solution of the
equation (16) is important to  the financial Black-Scholes model .

{\bf Acknowledgements.}\ \ We would like to thank our teacher
Processor Zhang Shiqing for his lectures on Functional Analysis,
meanwhile We would like to thank him for organizing the seminar on
financial mathematics and his many helpful discussions, suggestions
and corrections about this paper.


\vspace{0.3cm}
\begin{thebibliography}{99}
\bibitem{1} E. Zeidler, Applied Functional Analysics: Applications to Mathematical Physics, Spring-Verlag, New York-Berlin, 1995, P. 18-64.

\bibitem{2}Hui-Hsiung Kuo, Introduction to Stochastic
Integration, Springer, 2006, P. 43-61.

\bibitem{3} Bernt $\phi $ Ksendl,Stochastic Differential
Equations: An introduction with Applications, Springer-Verlag, 2005,
P.29.

\bibitem{4} Bharucha-Reid, A. P., Random Integral
Equations, Academic Press, New York, 1972.
\bibitem{5} Doob, J. Stochastic Processes, Wiley, New York, 1953.
\bibitem{6}Padgett, W. J. and Tsokos, C. P., On a stochastic integral of the Volterra type in telephone traffic
theory, J. Appl. Prob. 8, 1971, 269-275.

\bibitem{7}Tsokos, C. P. and Padgett, W. J., Random Integral
Equations with Application to Life Sciences and Engineering,
Academic Press, New York, 1974.

\bibitem{8}Rao, A. N. V. and Tsokos, C. P., On a class of
stochastic integral equation, Coll. Math. 35, 1976, 141-146.












\end{thebibliography}
\end{document}